\documentclass{article}
\usepackage[affil-it]{authblk}
\usepackage[utf8]{inputenc}
\usepackage[T1]{fontenc}
\usepackage[standard]{ntheorem}
\usepackage{amsmath}
\usepackage{calrsfs}

\usepackage{tikz}
\usepackage{thmtools}
\usepackage{thm-restate}
\usepackage{hyperref}
\usepackage{cleveref}
\usetikzlibrary{calc}

\tikzstyle{noeud}=[circle, inner sep=2, minimum size =3 pt, line width = 1pt, draw=black, fill=white]

\newcommand{\convexpath}[2]{
[   
    create hullnodes/.code={
        \global\edef\namelist{#1}
        \foreach [count=\counter] \nodename in \namelist {
            \global\edef\numberofnodes{\counter}
            \node at (\nodename) [draw=none,name=hullnode\counter] {};
        }
        \node at (hullnode\numberofnodes) [name=hullnode0,draw=none] {};
        \pgfmathtruncatemacro\lastnumber{\numberofnodes+1}
        \node at (hullnode1) [name=hullnode\lastnumber,draw=none] {};
    },
    create hullnodes
]
($(hullnode1)!#2!-90:(hullnode0)$)
\foreach [
    evaluate=\currentnode as \previousnode using \currentnode-1,
    evaluate=\currentnode as \nextnode using \currentnode+1
    ] \currentnode in {1,...,\numberofnodes} {
-- ($(hullnode\currentnode)!#2!-90:(hullnode\previousnode)$)
  let \p1 = ($(hullnode\currentnode)!#2!-90:(hullnode\previousnode) - (hullnode\currentnode)$),
    \n1 = {atan2(\y1,\x1)},
    \p2 = ($(hullnode\currentnode)!#2!90:(hullnode\nextnode) - (hullnode\currentnode)$),
    \n2 = {atan2(\y2,\x2)},
    \n{delta} = {-Mod(\n1-\n2,360)}
  in 
    {arc [start angle=\n1, delta angle=\n{delta}, radius=#2]}
}
-- cycle
}

\newtheorem{claim}[theorem]{Claim}

\newtheorem{problem}[theorem]{Problem}

\newtheorem{claimproof}{Proof of Claim}[claim]

\begin{document}

\title{Avoidance games are {\sf PSPACE}-Complete}

\author{Valentin Gledel
  \thanks{supported by the Kempe Foundation Grant No.~JCK-2022.1 (Sweden)}}
\affil{
Department of Mathematics
and Mathematical Statistics,\\
Umeå University, Sweden\\
\texttt{valentin.gledel@umu.se}}

\author{Nacim Oijid
  \thanks{supported by ANR-21-CE48-0001 project P-GASE}}
\affil{Univ. Lyon, Université Lyon 1,\\
LIRIS UMR CNRS 5205, F-69621, Lyon, France\\
\texttt{nacim.oijid@univ-lyon1.fr}}


\maketitle

\begin{abstract}
Avoidance games are games in which two players claim vertices of a hypergraph and try to avoid some structures. These games are studied since the introduction of the game of SIM in 1968, but only few complexity results are known on them. In 2001, Slany proved some partial results on Avoider-Avoider games complexity, and in 2017 Bonnet {\em et al.} proved that short Avoider-Enforcer games are Co-W[1]-hard. More recently, in 2022, Miltzow and Stojakovi{\'c} proved that these games are {\sf NP}-hard. As these games corresponds to the misère version of the well-known Maker-Breaker games, introduced in 1963 and proven {\sf PSPACE}-complete in 1978, one could expect these games to be {\sf PSPACE}-complete too, but the question remained open since then. We prove here that both Avoider-Avoider and Avoider-Enforcer conventions are {\sf PSPACE}-complete, and as a consequence of it that some particular Avoider-Enforcer games also are.

\end{abstract}

\section{Introduction}

\subsection{Related works}

Avoidance games belong to the class of positional games, that were introduced by Hales and Jewett in 1963~\cite{Hales1963} and popularised by Erd\H os and Selfridge in 1973~\cite{Erdos1973}. In this class of games, the board is a hypergraph and two players alternately select a vertex of the hypergraph that has not been claimed before. Winning conditions depend on the convention and are related to the hyperedges.
{\sc Tic-Tac-Toe} and {\sc Hex} are two famous examples of positional games. To learn more about positional games, we refer the reader to the recent survey of  Hefetz {\em et al.}~\cite{Hefetz2014}.

Among positional games, a natural dichotomy exists: on the one hand, there are games in which players seek to build a structure, and on the other hand, there are games in which players want to avoid a structure. The former set contains both Maker-Maker and Maker-Breaker conventions, in which the hyperedges are winning sets, and the player either want to fill a winning set (Maker role), or to play at least once in each of them (Breaker role). The latter contains Avoider-Avoider and Avoider-Enforcer conventions, that can be seen as the misère version of the former. In these games, the hyperedges are losing sets, and the players either want not to fill one losing set (Avoider role), or to force their opponent to fill one of them (Enforcer role).

When positional games were introduced, the focus was on Maker-Breaker games, {\em i.e.} games in which one player, Maker, aims to fill a hyperedge, and the second one, Breaker, wants to prevent it by claiming at least one vertex in each hyperedge. This convention is the most popular, and several games were studied according to this convention. In particular, the survey of Beck~\cite{Beck2008} presents several results obtained for Maker-Breaker games. The field of Maker-Breaker games is still well investigated today, and some Maker-Breaker games were introduced recently~\cite{Duchene2020, nenadov2016}.

The first Avoider-Avoider game was introduced in 1968 with the game of SIM and is presented in~\cite{simmons1968}, but the first study of the complexity of Avoidance games  was done by Schaefer in 1978~\cite{Schaefer1978}. Avoider-Enforcer games were introduced later by Lu in 1991~\cite{Lu1991,Lu1995} under the name of Antimaker-Antibreaker games and corresponds to the misère version of Maker-Breaker games. The standard name of this convention, Avoider-Enforcer, was popularised by Hefetz and different co-authors in 2007~\cite{hefetz2007-4, hefetz2007-2, hefetz2007-1, hefetz2007-3}. In this game, Enforcer wins if at some point during the game, Avoider has played all the vertices of a hyperedge, otherwise Avoider wins.

Even if most of the studies of positional games are focused on Maker-Breaker games, Avoider-Enforcer games have become more and more relevant: the famous Ramsey game were introduced in Avoider-Enforcer convention by Beck in 2002~\cite{beck2002} as a generalisation of SIM. As it was done in the Maker-Breaker convention, some games on graphs were introduced in Avoider-Enforcer or Avoider-Avoider conventions, where the loosing sets correspond to some structure in the graph, see~\cite{anuradha2008, barat2009, Grzesik2013, hefetz2010}.

In terms of complexity, an overview of the field is proposed by Demaine~\cite{demaine2001}. In positional games, as they are perfect information games, one player always has a winning strategy (or both players can ensure a draw). The natural decision problem related to games is thus: does the first player have a winning strategy ?
This problem was quickly proven to be {\sf PSPACE}-complete for Maker-Breaker games by Schaefer in 1978~\cite{Schaefer1978} even restricted to $11$-uniform hypergraphs (i.e. hypergraphs in which all hyperedges have size $11$). This bound was recently improved by Rahman and Watson in 2021~\cite{Rahman2021}, proving that the problem is still {\sf PSPACE}-complete if the hypergraph is $6$-uniform. These two proofs are very technical and a simpler proof of the {\sf PSPACE}-completeness was provided by Byskov in 2004~\cite{byskov2004}, proving in the same time that Maker-Maker games are also {\sf PSPACE}-complete. The complexity of Maker-Breaker games is still studied today, as Galliot {\em et al.}~\cite{galliot2022} have proven that the winner of a $3$-uniform Maker-Breaker game can be computed in polynomial time, but the gap between the complexity of $6$-uniform hypergraphs and $3$-uniform hypergraphs remains to be closed.

Despite the fact that Avoidance games were introduced at the same time as Maker-Breaker games, only partial results on complexity are known: determining the winner in Avoider-Avoider games, was proven to be {\sf PSPACE}-complete by Slany in 2002~\cite{slany2002} for endgames, i.e. games in which some vertices are already attributed to the players, but there is no result yet in the general case. Concerning Avoider-Enforcer,  Bonnet {\em et al.} in 2017~\cite{Bonnet2017} mentioned that the complexity of this problem is still open, when they proved that short games, i.e. games in which a player has only few moves to make, are co-W[1]-hard, with the number of moves taken as a parameter. The best known result today is due to Miltzow and Stojakovi{\'c} in 2022~\cite{Miltzow2022} that states the {\sf NP}-hardness of this decision problem and conjectures its {\sf PSPACE}-completeness.


\subsection{Presentation of the results}

The Avoider-Enforcer game is played as follows: given a hypergraph $H$, two players, called {\em Avoider} and {\em Enforcer}, alternately claim an unclaimed vertex of $H$ with Avoider starting. The game ends when all the vertices have been claimed. If Avoider has claimed all the vertices of a hyperedge, Enforcer wins. Otherwise, Avoider wins. The related decision problem is the following one.

\begin{problem}{\sc Avoider-Enforcer}

\noindent Input: A hypergraph $H$.

\smallskip

\noindent Output: True if and only if Avoider has a winning strategy in the Avoider-Enforcer game on $H$.
\end{problem}

This paper will focus on the proof of the following result:

\begin{restatable}{theorem}{maintheorem}\label{main theorem}
The {\sc Avoider-Enforcer} problem is {\sf PSPACE}-complete, even when the entry is restricted to hypergraphs with hyperedges of size at most $6$.
\end{restatable}

Our proof of \cref{main theorem} follows a similar idea to the proof of Rahman and Watson~\cite{Rahman2021} and the proof of Schaefer \cite{Schaefer1978}, by constructing some hyperedges forcing the order of the moves. Contrary to Maker-Breaker games, in Avoider-Enforcer convention, there is no vertex in which the players are urged to play, as in general, players do not want to move in avoidance games. The key idea of this reduction is to create some structures in which playing first is a losing move. In the provided construction, at any moment of the game, only few moves are not losing moves. Thus, we can control the vertices played by the two players.

The proof provided for {\sf PSPACE}-completeness of Avoider-Enforcer games, enables us to state the following corollary for Avoider-Avoider games that will also be proven later:

\begin{problem}{\sc Avoider-Avoider}

\noindent Input: A hypergraph $H$.

\smallskip

\noindent Output: True if and only if the second player has a winning strategy in the Avoider-Avoider game on $H$.
\end{problem}

\begin{restatable}{corollary}{avoideravoider}\label{Avoider-Avoider}
    The {\sc Avoider-Avoider} problem is {\sf PSPACE}-complete, even when the entry is restricted to $7$-uniform hypergraphs.
\end{restatable}


This paper is organised as follows. In section 2, we introduce two lemmas that will be used in the main proof of the results. In particular, we show that pairing strategies that are often used in Maker-Breaker conventions can also be applied to Avoider-Enforcer games.
Section 3 describes the reduction used to prove the {\sf PSPACE}-completeness and define an order on the move that we call the {\em legitimate order}. We also show in this section that the proof holds if both players follow the legitimate order.
In section 4, we show that if a player does not follow the legitimate order then it cannot be a disadvantage to the other player, completing the proof of the Theorem~\ref{main theorem}.
Finally, in section 5, we use Theorem~\ref{main theorem} to prove Corollary~\ref{Avoider-Avoider} and to show that the Avoider-Enforcer versions of the Domination game and of the $H$-game are {\sf PSPACE}-complete.

\section{Preliminaries}

In Maker-Breaker games, some moves appear to be better than others. For instance, if a vertex is in all the hyperedges, it is always an optimal move to play it for both players. We present here a similar result to prove that some moves are better than others. Next lemma was first introduced by Miltzow and Stojakovi\'c~\cite{Miltzow2022}, but we still provide a proof by consistency for the reader.

\begin{lemma} \label{lemma:included moves}
Let $H$ be a hypergraph, and $u,v$ two vertices of $H$ such that, for every hyperedge $e$ containing $u$, $e$ also contains $v$. If a player has a winning strategy, then this player has a winning strategy in which he never plays $v$ while $u$ is unclaimed.
\end{lemma}

\begin{proof}
Let $H$ be a hypergraphs and let $u,v$ be two vertices such that, for every hyperedge $e$ containing $u$, $e$ also contains $v$. Let $\mathcal{S}_1$ be a winning strategy for Avoider (Enforcer resp.). We define the strategy $\mathcal{S}_2$ as follows:

\begin{itemize}
    \item If $\mathcal{S}_1$ claims a vertex $w$ in $V\setminus \{v\}$, $\mathcal{S}_2$ claims $w$.
    \item If $\mathcal{S}_1$ claims $v$ while $u$ is unclaimed, claim $u$ and consider the strategy obtained if $v$ was claimed in $\mathcal{S}_1$. If at some point the other player plays $v$ continue as if he played $u$ in $\mathcal{S}_1$, and, if the strategy of $\mathcal{S}_1$ requires that you play $u$, play on $v$ and continue following $\mathcal{S}_1$.
\end{itemize}

Consider a sequence of moves played against $\mathcal{S}_2$. Consider the same sequence played against $\mathcal{S}_1$ by exchanging the roles of $u$ and $v$. Note $V_1$ the vertices claimed by Avoider (Enforcer resp.) in $\mathcal{S}_1$, and $V_2$ the vertices claimed by Avoider (Enforcer resp.) in $\mathcal{S}_2$ according to this sequence. Note that if $V_1 \neq V_2$, then $V_1 \setminus V_2 = \{v\}$ and $V_2 \setminus V_1 = \{u\}$.

As for every hyperedge $e$ such that $u \in e$, we have $v \in e$, in particular, we have for any hyperedge $e$, $|e \cap V_1| \ge |e \cap V_2|$. Thus, in this sequence, as $\mathcal{S}_1$ was a strategy for Avoider (Enforcer resp.) that do not claim all the vertices of a hyperedge (that do not claim a vertex in each hyperedge resp.), $\mathcal{S}_2$ neither. As this result holds for any sequence of moves, $\mathcal{S}_2$ is a winning strategy.
\end{proof}

The second tool we introduce here is {\it pairing strategies} in Avoider-Enforcer games. In Maker-Breaker convention, these strategies are often described by using the fact that a player can claim at least one vertex in each pair of vertices. Here in Avoider-Enforcer convention, the main idea of pairing strategies is that it is always possible to force the opponent to play at least once in each pair.

In this section, we will refer to the players by Alice and Bob, as the strategy can be applied both by Avoider and by Enforcer.

\begin{lemma}\label{pairing strategy}
Let $H = (V,E)$ be a hypergraph. Suppose that Alice plays last in $H$, i.e. if the game is played until all the vertices have been claimed, Alice will claim the last one.
Let $(a_1, b_1), \dots, (a_n, b_n)$ be pairwise disjoint pairs of vertices, and let $v \not \in \underset{i=1}{\overset{n}{\bigcup}} \{a_i,b_i\}$.

Alice has a strategy which forces Bob to play at least one vertex in each pair $(a_i,b_i)$. Bob has a strategy which forces Alice to play $v$ and at least one vertex in each pair $(a_i,b_i)$.
\end{lemma}

A strategy satisfying the hypothesis of Lemma~\ref{pairing strategy} will be called a {\em pairing strategy}.

\begin{proof}
We prove by induction that if Alice plays last on $H$ she can force Bob to play at least once in each pair of unplayed vertices $(a_i,b_i)$ :
\begin{itemize}
    \item If $|V| \in \{0,1\}$, then there can be no pair of distinct vertices $(a_i,b_i)$ and the strategy is trivially achieved.
    \item If $|V| = 2p+1$ with $p \geq 0$, then, since Alice plays last on $H$, it is Alice's turn to play and there is at least one vertex $x$ in $V$ that is in no pair and Alice can play it. Then, by induction hypothesis, Bob will play one vertex in each pair $(a_i, b_i)$ in $V\setminus \{x\}$.
    \item If $|V| = 2p$, with $p \geq 1$, then it is Bob's turn to play. If Bob plays a vertex $x$ that is in no pair $(a_i,b_i)$, then, by induction hypothesis, Alice can force Bob to play once in each pair $(a_i, b_i)$ in $V\setminus \{x\}$. If Bob plays on a vertex $a_i$ or $b_i$, then Alice can play in the other vertex of the pair. Bob played in the pair $(a_i, b_i)$ and by induction hypothesis Alice can force Bob to play once in each pair $(a_j, b_j)_{j \neq i}$ in $V \setminus \{a_i,b_i\}$.
\end{itemize}

We prove by induction that if Alice plays last on $H$, Bob can force Alice to play at least once in each pair of unplayed vertices $(a_i,b_i)$, and to play on the vertex $v$ :
\begin{itemize}
    \item If $|V| = 0$, the situation cannot exist as there can be no vertex $v$.
    \item If $|V| = 1$, it must be that $V = \{v\}$ and that there is no pair of vertices $(a_i,b_i)$ and since it is the last move, it is Alice's turn, and she has to play in $v$. Thus, she played on $v$ and in each pair.
    \item If $|V| = 2p$ with $p \geq 1$, then, since Alice plays last on $H$, it is Bob's turn to play and there is at least one vertex $x$ in $V\setminus \{v\}$ that is in no pair and Bob can play it. Then, by induction hypothesis, Alice will play on $v$ and on one vertex in each pair $(a_i, b_i)$ in $V\setminus \{x\}$.
    \item If $|V| = 2p+1$, with $p \geq 1$, then it is Alice's turn to play. If Alice plays a vertex $x$ that is neither $v$ neither in a pair $(a_i,b_i)$, then by induction hypothesis Bob can force Alice to play one $v$ and at lea  st once in each pair $(a_i, b_i)$ in $V\setminus \{x\}$. If Alice plays on a vertex $a_i$ or $b_i$ then Bob can play in the other vertex of the pair. Alice played in the pair $(a_i, b_i)$ and by induction hypothesis Bob can force Alice to play on $v$ and once in each pair $(a_j, b_j)_{j \neq i}$ in $V \setminus \{a_i,b_i\}$. If Alice plays on $v$ then either there is no pair $(a_i,b_i)$ and then Alice played on $v$ and in each pair, or there is at least one pair $(a_i,b_i)$. Bob can then play on $a_i$, by induction hypothesis Bob can force Alice to play on $b_i$ and in each pair $(a_j, b_j)_{j \neq i}$. Since Alice played on $v$, and will be forced to play on $b_i$ and in each pair $(a_j, b_j)_{j \neq i}$ the result holds.
\end{itemize}
\end{proof}

\section{Proof of the main theorem}

In this section, we begin the proof of Theorem~\ref{main theorem} by describing the reduction from $3$-QBF, introducing an order on the move of Avoider and Enforcer called the {\em legitimate order} and giving a sketch of the general proof.

\maintheorem*

The first step of the proof is to prove that this game is in {\sf PSPACE}.

\begin{lemma}\label{inPspapce}
The {\sc Avoider-Enforcer} problem is in {\sf PSPACE}.
\end{lemma}

\begin{proof}
Let $H = (V,E)$ be a hypergraph. As the players are not allowed to play an already claimed vertex, any game ends after at most $|V|$ moves. Therefore, according to Lemma~2.2 of Schaefer~\cite{Schaefer1978}, as the game has a polynomial length and a polynomial number of moves, its winner can be computed with polynomial space.
\end{proof}

\subsection{Construction of the hypergraph}\label{subsection construction}

We reduce the problem {\sc 3-QBF} to an {\sc Avoider-Enforcer} game. This problem has been proven {\sf PSPACE}-complete by Stockemeyer and Meyer \cite{stockmeyer1973}, and we use the gaming version of this problem as it was formulated by Rahman and Watson~\cite{Rahman2021}. The game is played on a quantified formula $\varphi$ of the form $\forall X_1 \exists X_2 \dots \forall X_{2n-1} \exists X_{2n} \psi$, with $\psi$ a $3$-SAT formula. Alternately, two players, namely Falsifier and Satisfier, chose valuation for the variables, Falsifier for the odd variables (quantified with a $\forall$) and Satisfier for the even ones (quantified with a $\exists$). When all the variables have a valuation Satisfier wins if $\psi$ is satisfied, otherwise, Falsifier win.

\begin{problem}{3-QBF}

\noindent Input: A $3$-SAT quantified formula $\varphi$ of the form $\forall X_1 \exists X_2 \dots \forall X_{2n-1} \exists X_{2n} \psi$.

\smallskip

\noindent Output: True if and only if Satisfier has a winning strategy in the 3-QBF game on $\phi$
\end{problem}

A round in a $3$-QBF formula corresponds to the steps $i$ during which Falsifier gives a valuation to $X_{2i-1}$ and then Satisfier gives a valuation to $X_{2i}$. In this reduction, any round corresponds to ten vertices and eight hyperedges. Four of the ten vertices are $\{x_{2i-1}, \overline{x_{2i-1}}, x_{2i}, \overline{x_{2i}}\}$, and the six others are $u_{6i-5}, u_{6i-4}, u_{6i-3}, u_{6i-2}, u_{6i-1}, u_{6i}$. The eight hyperedges are constructed as follows:

\vspace{-0.4cm}

\begin{align*}
    A_{2i} &=  \hspace{.1cm} ( x_{2i}, \overline{x_{2i}}, u_{6i+1}, u_{6i+3} ) &\\
    C^+_{6i} &=  \hspace{.1cm} ( u_{6i}, u_{6_i+1}, u_{6i+3}, x_{2i} ) &\\
    C^-_{6i} &=  \hspace{.1cm} ( u_{6i}, u_{6_i+1}, u_{6i+3}, \overline{x_{2i}} ) &\\
    C^+_{6i-2} &=  \hspace{.1cm} ( u_{6i-2}, u_{6_i-1}, u_{6i+1}, x_{2i} ) &\\
    C^-_{6i-2} &=  \hspace{.1cm} ( u_{6i-2}, u_{6_i-1}, u_{6i+1}, \overline{x_{2i}} ) &\\
    B_{2i-1} &=  \hspace{.1cm} ( x_{2i-1}, \overline{x_{2i-1}}, u_{6i-1} ) &\\
    C^+_{6i-4} &=  \hspace{.1cm} ( u_{6i-4}, u_{6_i-3}, u_{6i-1}, x_{2i-1} ) &\\
    C^-_{6i-4} &=  \hspace{.1cm} ( u_{6i-4}, u_{6_i-3}, u_{6i-1}, \overline{x_{2i-1}} ) & 
\end{align*}

\medskip

\noindent Moreover, for each clause $F_j = l_{j_1} \vee l_{j_2} \vee l_{j_3} \in \psi$ where the vertices $l_{j_1}, l_{j_2}$ and $l_{j_3}$ are literals either positive or negative, we add a hyperedge
$D_j$ containing six vertices. For $k = 1,2,3$, if $l_{j_k}$ is a positive variable $X_i$, then $x_i$ is in $D_j$, if $l_{j_k}$ is a negative one $\neg X_i$, then $\overline{x_i}$ is in $D_j$. If $j_k$ is odd, then $u_{6j_k -1}$ is in $D_j$, if $j_k$ is even, then $u_{6 j_k +1}$ is in $D_j$.



\medskip

\noindent Finally, the CNF game $\phi$ is reduced to the hypergraph $H = (V,E)$ with

\vspace{-0.25cm}

$$    V =  \hspace{.1cm} \left \{ \{x_i\}_{1 \le i \le 2n} \cup \{ \overline{x_i}\}_{1 \le i \le 2n} \cup \{ u_j \}_{1 \le j \le 6n} \right \}$$ $$  E =  \hspace{.1cm}  \left \{ \{A_{2i}\}_{1 \le i \le n} \cup \{C^+_{2i}\}_{1 \le i \le 3n} \cup \{C^-_{2i}\}_{1 \le i \le 3n} \cup \{B_{2i-1}\}_{1 \le i \le n} \cup \{D_j\}_{1 \le j \le m} \right \} 
$$

With this construction, we say that Avoider and Enforcer follow a {\em legitimate order} if during each round $i$ the moves are made in the following way, with each round played in increasing order:

\noindent {\bf Legitimate order during round $i$}
\begin{enumerate}
  \setlength\itemsep{-.25em}
    \item Avoider starts and plays $u_{6i-5}$. 
    \item Enforcer plays $u_{6i-4}$.
    \item Avoider plays $u_{6i-3}$. 
    \item Enforcer plays one of $x_{2i-1}$ or $\overline{x_{2i-1}}$.
    \item Avoider plays the remaining vertex in $(x_{2i-1}, \overline{x_{2i-1}})$
    \item Enforcer plays $u_{6i-2}$. 
    \item Avoider plays $u_{6i-1}$.
    \item Enforcer plays $u_{6i}$. 
    \item Avoider plays one of $x_{2i}$ or $\overline{x_{2i}}$.
    \item Enforcer plays the remaining vertex in $(x_{2i}, \overline{x_{2i}})$.
\end{enumerate}

\subsection{Sketch of the proof}

To prove that, with our construction, Avoider wins the {\sc Avoider-Enforcer} game, if and only if Satisfier wins the QBF game, we first prove that this statement is true if the order of the move is legitimate, as the moves will give the valuation obtained in QBF. To force the players to play in the legitimate order, the main idea of the construction is that players want to play some vertices as late as possible. Therefore, we prove that it is always optimal to respect the legitimate order of the moves. We introduce the three following lemmas that will be proved in the next section.

\begin{restatable}{lemma}{lemmaorder}
\label{lem:order_resp}
When the game is restricted to the legitimate order, Avoider has a winning strategy in the Avoider-Enforcer game on $H$ if and only if Satisfier has a winning strategy for the {\sc 3-QBF} game on $\varphi$.
\end{restatable}

\begin{restatable}{lemma}{falsifierwin}
\label{lem:false_win}
If Enforcer has a winning strategy in $\varphi$ when the legitimate order is respected by the two players, then he has a winning strategy in $H$.
\end{restatable}

\begin{restatable}{lemma}{satisfierwin}
\label{lem: satis win}
If Avoider has a winning strategy in $\varphi$ when the legitimate order is respected by the two players, then she has a winning strategy in $H$.
\end{restatable}

We first admit these lemmas and we prove \cref{main theorem}.

\begin{proof}

Let $\varphi$ be a boolean formula. Consider the hypergraph $H$ obtained from $\varphi$ by following the construction of \cref{subsection construction}. This construction has polynomial size, and the hypergraph $H$ has all its hyperedges of size at most $6$. According to \cref{lem:order_resp}, when the order is respected, if Satisfier (Falsifier resp.) has a winning strategy in $\varphi$, Avoider (Enforcer resp.) has a winning strategy in $H$. Thus, according to \cref{lem: satis win} ( \cref{lem:false_win} resp.), if Avoider (Enforcer resp.) has a winning strategy on $H$ when the legitimate order is respected, she (he resp.)  has one in general in $H$. Thus, Satisfier wins on $\varphi$ if and only if, Avoider wins on $H$. Therefore, the {\sc Avoider-Enforcer problem} is {\sf PSPACE}-complete.
\end{proof}

\subsection{Game in legitimate order} \label{order respected}

In this section, we suppose that both players follow a legitimate order of moves.

If the order of moves is legitimate, the only choices available for Avoider and Enforcer are on the vertices $x_i$ and $\overline{x_i}$. For each $1 \le i \le 2n$, Avoider plays one of $x_i, \overline{x_i}$ and Enforcer the other. Therefore, if both Avoider and Enforcer plays one vertex in $\{x_i, \overline{x_i}\}$, we define the {\em underlying valuation} given to $\psi$ as the following one:

\begin{align*}
    X_i = \left\{
    \begin{array}{ll}
        True & \mbox{if Avoider has played $\overline{x_i}$ and Enforcer has played $x_i$} \\
        False & \mbox{if Avoider has played $x_i$ and Enforcer has played $\overline{x_i}$}
    \end{array}
\right. 
\end{align*}

\lemmaorder*

\begin{proof}

Consider a game played on $H$ for which both Avoider and Enforcer respected the legitimate order through the whole game. 

\begin{claim}
Avoider won the game on $H$ if and only if the formula $\psi$ is satisfied by the underlying valuation of the $X_i$s.
\end{claim}

\begin{claimproof}
Since the legitimate order is respected, Enforcer played all the vertices $u_{2i}$ and thus played at least once in all the hyperedges $C^+_{2i}$ and $C^-_{2i}$. Moreover, for each pair of variables $(x_i,\overline{x_i})$, Enforcer played one of the vertices of the pair and so played at least one vertex in all the hyperedges $A_i$ and $B_i$.
Thus, the only hyperedges that could possibly be fully played by Avoider are the hyperedges $D_j$.

Since, in the legitimate order, Avoider plays on all the vertices $u_{2i+1}$, a hyperedge $D_j$ is fully played by Avoider if and only if she played on all the vertices $x(l_k)$ for $l_k \in F_j$. If this is the case, then this means that the formula $\psi$ is not satisfied by the underlying valuation because the clause $F_j$ has all its literals assigned to False. On the contrary, if the formula $\psi$ is satisfied by the underlying valuation, then, for all clause $F_j$, at least one of its literal is assigned to True and so Enforcer played at least once in each hyperedge $D_j$.

Therefore, Avoider won the game on $H$ if and only if $\psi$ is satisfied.
\end{claimproof}

Suppose Satisfier has a winning strategy $\mathcal{S}$ on $\varphi$. We define a strategy for Avoider as follows: Whenever Avoider has to play a vertex $x_{2k}$ or $\overline{x_{2k}}$, Avoider considers the underlying valuation given to the $X_i$s with $i < 2k$. Then, if Satisfier had put $X_{2k}$ to True, she plays $\overline{x_{2k}}$. Otherwise, she plays $x_{2k}$. With this strategy, at the end of the game, the underlying valuation of the variables of $H$ will be equivalent to the valuation given by the game that Satisfier played on $\varphi$. Since Satisfier has a winning strategy on $\varphi$, the underlying valuation satisfies $\psi$ and so Avoider win the game.

Similarly, if Falsifier has a winning strategy, Enforcer can follow the strategy in such a way that at the end the valuation of variables in the game played by Falsifier correspond to the underlying valuation in $H$. Since Falsifier has a winning strategy, Enforcer wins the game on $H$.
\end{proof}

\section{Proofs of Lemma~\ref{lem:false_win} and Lemma~\ref{lem: satis win}}

The first part of our constructions showed that, if the legitimate order is respected, Avoider wins if and only if $\varphi$ is satisfied by the underlying valuation provided by the moves. We now prove that, if a player has a winning strategy, he has one that respects the legitimate order of the moves. We introduce here sets of variables in which players do not want to play. These sets will be the main point of the proofs of \cref{lem:false_win} and \cref{lem: satis win}.

For $i = 1$ to $4n$, 
we define the set of vertices $S_i$ as $S_{4n}  = \{u_{6n},x_{2n},\overline{x_{2n}}\}$ and for $i < 4n$:
\begin{itemize}
    \item if $i=4k$, $S_i = \{u_{6k},x_{2k},\overline{x_{2k}},u_{6k+1}\} \cup S_{i+1} $
    
    \item if $i=4k-1$, $S_i = \{u_{6k-2},u_{6k-1}\} \cup S_{i+1}$
    
    \item if $i=4k-2$, $S_i = \{x_{2k-1},\overline{x_{2k-1}}\} \cup S_{i+1}$
    
    \item if $i=4k-3$, $S_i = \{u_{6k-4},u_{6k-3}\} \cup S_{i+1}$
    
\end{itemize}


\subsection{Proof of Lemma~\ref{lem:false_win}}

\falsifierwin*

\begin{proof}

Suppose Falsifier has a winning strategy in $\varphi$. Consider a strategy for Enforcer in which he plays according to the legitimate order until Avoider does not. If Avoider respects the order until all the vertices are played, by Lemma~\ref{lem:order_resp}, Enforcer wins. Otherwise, the proof of the following claim provides a winning strategy to Enforcer.

\begin{claim}
If, during the game, Avoider plays in a set $S_i$ in which Enforcer has not played yet, then, after this move, Enforcer has a strategy to win the game.
\end{claim}

\begin{claimproof}
The proof is by induction on $i$. 

First, notice that each $S_i$ has an odd number of vertices and, as the total number of vertices in $H$ is $10n$, there is also an odd number of vertices outside $S_i$. Therefore, if Avoider plays first in an $S_i$, Enforcer answers by playing an arbitrary vertex that is not in $S_i$ and considers an arbitrary pairing outside $S_i$, which exists as there is an even number of vertices outside $S_i$ after his move. This way, Avoider has to be the next player to play in $S_i$.

{\bf Base cases:}
\begin{itemize}
    \item Case $i= 4n$: If Avoider plays first is $S_{4n}$, by pairing the two other vertices in $S_{4n}$, by using Lemma~\ref{pairing strategy}, Enforcer can force Avoider to play another vertex in $S_{4n}$. Hence, as $(u_{6n},x_{2n})$, $(u_{6n},\overline{x_{2n}})$ and $(x_{2n}, \overline{x_{2n}})$ are three hyperedges, Avoider will claim the two vertices of one of them and then lose.
    \item Case $i=4n-1$: As shown previously, Enforcer has as strategy such that Avoider is the next player to play in $S_{4n-1}$. If Avoider has played at least one of its two first moves in $S_{4n}$, she has lost by the case $i=4n$. Otherwise, she has played exactly $u_{6n-2}$ and $u_{6n-1}$. In this case, Enforcer plays on $u_{6n}$ and pairs $x_{2n}$ and $\overline{x_{2n}}$ and by Lemma~\ref{pairing strategy} he forces Avoider to claim all the vertices of $C^+_{6n-2}$ or $C^-_{6n-2}$.
\end{itemize}

{\bf Induction steps:}
Suppose that the first time Avoider does not respect the order of the move, she plays in a set $S_i$ for $i \le 4n-2$. If the second move of Avoider in $S_i$ is in $S_{i+1}$, Enforcer wins by induction hypothesis. Thus, we can suppose that Avoider has played two vertices in $S_i \setminus S_{i+1}$. Moreover, as Enforcer has arbitrarily paired the vertices outside $S_i$, we describe here the strategy in $S_i$, and Enforcer plays according to the pairing outside $S_i$. This strategy ensure that the moves in $S_i$ alternate between both players.

\begin{itemize}
    \item Case $i= 4k$: Avoider has played twice in $\{u_{6k}, x_{2k}, \overline{x_{2k}}, u_{6k+1}\}$. At least one of $\{u_{6k}, x_{2k}, \overline{x_{2k}}\}$ is available. Enforcer plays it. Avoider has to play a third vertex in this quadruple, otherwise, she plays first in $S_{i+1}$ and loses by induction, and necessarily one of the three vertices she has played is $u_{6k+1}$. Enforcer plays $u_{6k+2}$. Avoider either plays first in $S_{i+2}$ and loses by induction hypothesis, or plays $u_{6k+3}$. At this moment, Avoider has played on the vertices $u_{6k+1}$ and $u_{6k+3}$, and two of the vertices of $\{u_{6k}, x_{2k}, \overline{x_{2k}}\}$. So she has completed one of the hyperedge $C^+_{6k} = (u_{6k}, u_{6k+1}, u_{6k+3}, x_{2k})$, $C^-_{6k} = (u_{6k}, u_{6k+1}, u_{6k+3}, \overline{x_{2k}})$ or $A_{2k} = (x_{2k}, \overline{x_{2k}}, u_{6k+1}, u_{6k+3}) $.
    \item Case $i=4k-1$: Avoider has played $u_{6k-2}$ and $u_{6k-1}$. Enforcer plays $u_{6k}$. Avoider has to play on vertex in $\{x_{2k}, \overline{x_{2k}}, u_{6k+1}\}$, otherwise she plays first in $S_{i+2}$ and loses by induction. Enforcer plays either $x_{2k}$ or $\overline{x_{2k}}$, as at least one of them is available. If Avoider plays a vertex in $S_{i+2}$ she loses by induction. So she has to play the last vertex available in $S_i \setminus S_{i+1}$. With this strategy, Avoider has necessarily played $u_{6k+1}$ and one of $x_{2k}$ and $\overline{x_{2k}}$. Thus, she has played all the vertices of either $C^+_{6k-2} = (u_{6k-2}, u_{6k-1}, u_{6k+1}, x_{2k})$ or $C^-_{6k-2} = (u_{6k-2}, u_{6k-1}, u_{6k+1}, \overline{x_{2k}})$.
    \item Case $i=4k-2$: Avoider has played $x_{2k-1}$ and $\overline{x_{2k-1}}$. Enforcer plays $u_{6k-2}$. Either Avoider plays first in $S_{i+2}$ and loses by induction, or she plays $u_{6k-1}$, the last available vertex in $S_{i+1}$ and loses by having played all the vertices in $B_{2k-1} = (x_{2k-1}, \overline{x_{2k-1}}, u_{6k-1})$.
    \item Case $i=4k-3$: Avoider has played $u_{6k-4}$ and $u_{6k-3}$. Enforcer plays $x_{2k-1}$. Avoider has to play $\overline{x_{2k-1}}$, otherwise she plays first in $S_{i+2}$ and loses by induction. Then Enforcer plays $u_{6k-2}$. If Avoider plays in $S_{i+3}$ she loses by induction. The last vertex available in $S_i \setminus S_{i+3}$ is $u_{6k-1}$, and if Avoider plays it, she loses by playing all the vertices $C^-_{6k-4} = (u_{6k-4}, u_{6k-3}, u_{6k-1}, \overline{x_{2k-1}})$. 
\end{itemize}

By applying this induction, at any moment of the game, if Avoider plays first in a set $S_i$, she looses.
\end{claimproof}

Finally, if Enforcer has played according to the legitimate order, at any moment of the game, Avoider has to play in a set $S_i$ in which Avoider has already played. Therefore, she has to respect the order of the moves. The only moment when she can change this order is by playing $u_{6k+1}$ instead of one of the vertices $x_{2k}, \overline{x_{2k}}$. But if she does so, Enforcer can play one of them, for instance $x_{2k}$, and Avoider will be forced to play $\overline{x_{2k}}$. If this happens, everything happens as if Avoider has played $\overline{x_{2k}}$ first and $u_{6k+1}$ after. The strategy can then continue as if the order has been respected, as this does not change the order in which the valuation of the variables is chosen in $\varphi$, and Falsifier had a winning strategy if Satisfier had put $X_{2k}$ to True by hypothesis.

To conclude, if Falsifier has a winning strategy in $\varphi$, Enforcer has a winning strategy in $H$.
\end{proof}

\subsection{Proof of Lemma~\ref{lem: satis win}}

In this section, we prove that if Satisfier has a winning strategy in $\varphi$, Avoider has one in $H$, even if Enforcer does not respect the order. The main idea of the strategy is to respect the order, and if Enforcer does not respect the order, Avoider has a pairing strategy to force enforcer to claim some vertices $x_i$, or the odd $u_j$ that follows it. By construction, any hyperedge containing a vertex $x_i$ contains also the next odd vertex $u_i$ in the order, and this will prove that when Enforcer do not respect the order, it will benefit Avoider.

\satisfierwin*

\begin{proof}
    
If Satisfier has a winning strategy $\mathcal{S}$ in $\varphi$, and if Enforcer respects the order, Avoider has a winning strategy according to Lemma~\ref{lem:order_resp}.

While Enforcer respects the legitimate order, Avoider also respects it. Suppose that at any moment of the game, Enforcer does not respect the legitimate order. Denote by $y_A$ the vertex he would have played according to the legitimate order, and by $y_E$ the vertex he has played instead. First, note that, according to Lemma~\ref{lemma:included moves}, we can suppose that $y_E$ is not a vertex $u_j$ with $j$ odd. Indeed, for each vertex $u_{2i+1}$, the hyperedges that contain the previous vertex in the legitimate order (if this vertex is a vertex $x_j$ or $\overline{x_j}$ this is true for any of them), also contain $u_{2i+1}$. If $y_A$ is a vertex $x_{2i}$ or $\overline{x_{2i}}$, Avoider pairs it with $u_{6i+1}$ and continues as if Enforcer should have played $u_{6i+2}$. As any hyperedge containing $x_{2i}$ or $\overline{x_{2i}}$, also contains $u_{6i+1}$ which is the next vertex Avoider should have played according to the legitimate order, it benefits her if Enforcer finally plays $u_{6i+1}$ according to Lemma~\ref{lemma:included moves}.

Denote by $k$ the smallest integer such that $y_E \notin S_k$, and by $k'$ the largest integer such that $y_A \in S_{k'}$. Note that all the vertices outside $S_{k'}$ have already been played or are paired, and that $S_{k'} \setminus S_k$ is then the set of vertices perturbed by the move of Enforcer. We consider $k = 4n+1$ if $y_E \in S_{4n}$. As all the set $S_i$ have an odd number of variables, Avoider knows that the number of remaining vertices outside $S_k$ is odd. (as an even number of moves have been done in it). According to $\mathcal{S}$, Satisfier has a winning strategy starting with the value already given to the variables $x_i$  and $\overline{x_i}$ that are outside $S_{k'}$. Avoider considers than arbitrary moves for Falsifier in $\mathcal{S}$ and the corresponding moves for Satisfier until all the vertices $x_i$s and $\overline{x_i}$s in $S_k$ are played in $\mathcal{S}$. According to these moves, we will denote by $x^E_j$ the vertex among $(x_j, \overline{x_j})$ played by Enforcer and by $x^A_j$ the vertex played by Avoider, such that their underlying valuation is the one desired by $\mathcal{S}$. 

Avoider will then play a strategy different in $S_k$ and outside $S_k$:
\begin{itemize}
    \item In $H\setminus S_{k}$, Avoider plays $y_A$, the vertex that Enforcer should have played according to the order, then plays according to a pairing strategy, that is presented in the next paragraph.
    \item In $S_k$, Avoider considers the strategy she would have played if all the vertices outside $S_k$ were played according to the legitimate order, with the valuation she considered in $\mathcal{S}$.
\end{itemize}

The pairing we define is the following one: $(u_{6i-4}, x_{2i-1}^A)$, $(u_{6i-3}, u_{6i-6})$, $(x_{2i-1}^E, u_{6i-1})$, $(u_{6i-2}, x_{2i}^A)$, $(u_{6i+1}, x_{2i}^E)$. This pairing starts at the vertex $y_A$ and we consider only pairs containing a vertex outside $S_k$. Note that, by construction, exactly one vertex of this pairing is already played, and exactly one is in $S_k$. Therefore, to make the pairing contain only vertices not played and outside $S_k$, some modifications are done. These modifications are presented in Figure~\ref{fig:changes matching}. By applying Lemma~\ref{pairing strategy}, Avoider can ensure that Enforcer plays at least one in each of these pairs.

\begin{figure}
    \centering
\begin{tabular}{ |c|c| } \hline
 $y_A$  & changes \\ 
 \hline
  $u_{6i-4}$ & $u_{6i-3} \longleftrightarrow x_{2i-1}^A$ \\   
  $x_{2i-1}$ or $\overline{x_{2i-1}}$ & $x^*_{2i-1} \longleftrightarrow u_{6i-1}$ \\   
 $u_{6i-2}$ & $u_{6i-1} \longleftrightarrow x_{2i}^A$ \\   
  $u_{6i}$ & $u_{6i+3} \longleftrightarrow x_{2i}^A$ \\ \hline
\end{tabular}
\begin{tabular}{ |c|c| } \hline
 $y_E$&  changes \\ \hline
 $u_{6i-4}$ & no changes \\  
 $x_{2i-1}$ or $\overline{x_{2i-1}}$ & $x_{2i-1}^*  \longleftrightarrow u_{6i-4}$  \\   
 $u_{6i-2}$  & no changes \\    
 $u_{6i}$  & no changes\\  
 $x_{2i}$ or $\overline{x_{2i}}$ & $x_{2i}^*  \longleftrightarrow u_{6i+1}$, $u_{6i-2} \longleftrightarrow u_{6i}$  \\\hline    
\end{tabular} \hfill
    
    \caption{Changes of the matching. ($x_k^*$ refers to the variable between $x_k$ and $\overline{x_k}$ that has not been played ) }
    \label{fig:changes matching}
\end{figure}

\begin{claim}\label{strat ABC}
    The pairing strategy ensures that Enforcer plays at least once in each hyperedge $A_i$, $B_i$ or $C_i$ containing at all their vertices in $S_{k'}$ and at least one outside $S_k$.
\end{claim}

\begin{claimproof}

First, if the hyperedge contain no vertex whose pairing has been modified because of their appurtenance to $y_A$ or $y_E$, it contains two paired vertices. We show in bold text the paired vertices:

\vspace{-.4cm}

\begin{align*}\label{pairing Avoider}
    A_{2i} &=  \hspace{.1cm} ( x^A_{2i}, {\bf x^E_{2i}}, {\bf u_{6i+1} }, u_{6i+3} ) &\\
    C^A_{6i} &=  \hspace{.1cm} ( {\bf u_{6i}}, u_{6_i+1}, {\bf u_{6i+3}}, x^A_{2i} ) &\\
    C^E_{6i} &=  \hspace{.1cm} ( u_{6i}, {\bf u_{6_i+1}}, u_{6i+3}, {\bf x^E_{2i}} ) &\\
    C^A_{6i-2} &=  \hspace{.1cm} ( {\bf u_{6i-2}}, u_{6_i-1}, u_{6i+1}, {\bf x^A_{2i}} ) &\\
    C^E_{6i-2} &=  \hspace{.1cm} ( u_{6i-2}, u_{6_i-1}, {\bf u_{6i+1}}, {\bf x^E_{2i}} ) &\\
    B_{2i-1} &=  \hspace{.1cm} ( x^A_{2i-1}, {\bf x^E_{2i-1}}, {\bf u_{6i-1}} ) &\\
    C^A_{6i-4} &=  \hspace{.1cm} ( {\bf u_{6i-4}}, u_{6_i-3}, u_{6i-1}, {\bf x^A_{2i-1}} ) &\\
    C^E_{6i-4} &=  \hspace{.1cm} ( u_{6i-4}, u_{6_i-3}, {\bf u_{6i-1}}, {\bf x^E_{2i-1}} ) & 
\end{align*}

\smallskip

For the first hyperedges of the matching, there are two paired vertices:
\begin{itemize}
    \item If $y_A = u_{6i-4}$, only the hyperedges $C^E_{6i-4}$ and $C^A_{6i-4}$ are concerned by the changes. In the former $x^E_{2i-1}$ is paired with $u_{6i-1}$, in the latter $x^A_{2i-1}$ is paired with $u_{6i-3}$. 
    \item If $y_A\in \{ x_{2i-1}, \overline{x_{2i-1} } \}$, the only hyperedges concerned by the change is $B_{2i-1}$. In it, the other vertex in $ \{ x_{2i-1}, \overline{x_{2i-1} } \}$ is paired with $u_{6i-1}$
    \item If $y_A = u_{6i-2}$, only the hyperedges $C^E_{6i-2}$ and $C^A_{6i-2}$ are concerned by the changes. In the former $x^E_{2i}$ is paired with $u_{6i+1}$, in the latter $x^A_{2i}$ is paired with $u_{6i-1}$. 
    \item If $y_A = u_{6i}$, only the hyperedges $C^E_{6i}$ and $C^A_{6i}$ are concerned by the changes. In the former $x^E_{2i}$ is paired with $u_{6i+1}$, in the latter $x^A_{2i}$ is paired with $u_{6i+3}$. 
\end{itemize}

For the last hyperedges that contains vertices of the matching, the following happens:

\begin{itemize}
    \item If $y_E = u_{6i-4}$, the pairing stops at $u_{6i-3}$. The only two hyperedges that contains at least one vertex in $S_k$ and one vertex outside $S_k$ are $C^+_{6i-4}$ and $C^-_{6i-4}$, in which Enforcer has played $y_E$.
    \item If $y_E = x_{2i -1}$ or $\overline{x_{2i-1}}$, the pairing stops after the second vertex in $\{x_{2i-1}, \overline{x_{2i-1}}\}$. The only one hyperedge containing at least one vertex in $S_k$ and one outside $S_k$ is $B_{2i-1}$, in which Enforcer has already played $y_E$.
    \item If $y_E = u_{6i-2}$, the pairing stops at $u_{6i-1}$. The only two hyperedges that contains at least one vertex in $S_k$ and one vertex outside $S_k$ are $C^+_{6i-2}$ and $C^-_{6i-2}$, in which Enforcer has played $y_E$.
    \item If $y_E = u_{6i}$, the pairing stops at $u_{6i+1}$. The three hyperedges that contains both vertices in $S_k$ and vertices outside $S_k$ are $C^+_{6i}$, $C^-_{6i}$ and $A_{2i}$. In  $C^+_{6i}$, $C^-_{6i}$, Enforcer has played $y_E$, and in $A_{2i}$, Enforcer will play one of $x^E_{2i}$ or $u_{6i+1}$ as these two vertices are paired together.
    \item If $y_E = x_{2i}$ or $\overline{x_{2i}}$, the pairing stops at $u_{6i+1}$. The three hyperedges that contains vertices inside $S_k$ and outside $S_k$ are $C^+_{6i}$, $C^-_{6i}$ and $A_{2i}$. As the second vertex in $\{x_{2i}, \overline{x_{2i}}\}$ is paired with $u_{6i+1}$, either Enforcer has played both $x_{2i}$ and $\overline{x_{2i}}$, and any of these three hyperedges contain at least one of them; or Enforcer has played $u_{6i+1}$ which is in these three hyperedges.
\end{itemize}

If the pairing stops because it goes until the end (i.e. $k = 4n +1$), one vertex is not paired. According to Lemma~\ref{pairing strategy}, as Enforcer plays the last move in $H$, Avoider can force him to play it and still play once in each pair of the pairing.

Finally, in any hyperedge $A_i, B_i$ or $C_i$ containing at least one vertex of the matching, Enforcer has played at least one vertex. 
\end{claimproof}

Now, we can prove that the strategy we defined for Avoider is a winning strategy. In all the hyperedges $A_i, B_i$ or $C_i$, Enforcer played at least once: If Enforcer has respected the order until vertices of this hyperedge are reached, he has to play in it, otherwise by Claim~\ref{strat ABC}, Avoider can force Enforcer to play in it as this hyperedge is considered in a set of hyperedges in which Enforcer has not respected the order, as the only hyperedges in which Enforcer is not the first player to claim a vertex according to the order are the $A_{2i}$, and in them, the vertex $x_{2i}^E$ is always paired with $u_{6i+1}$.

In the hyperedges $D_j$, as the strategy defined by $\mathcal{S}$ is a winning strategy, Satisfier has a strategy to force at least one literal $l_i$ to be True in $F_j$. By construction, if when the vertex $x(l_i)$ has to be played, the order was respected, Enforcer has played it. If it has not, then Avoider has paired it with the vertex $u(l_i)$. In both cases, Enforcer has played in $D_j$.
\end{proof}

\section{Applications}

\subsection{From 6-hypergraphs to 6-uniform hypergraphs}

The construction provided in Section~\ref{subsection construction} provided a hypergraphs in which all hyperedges have size at most six. We prove here that we can suppose that all hyperedges have size six without changing the outcome.

A hypergraph $H = (V,E)$ is a {\it $k$-hypergraph} if each edge $e\in E$ has size at most $k$. It is said to be {\it $k$-uniform} if each edge $e \in E$ has size exactly $k$.

\begin{lemma}\label{to6uniform}
Let $H = (V,E)$ be a $k$-hypergraph. Let $m = \underset{e \in E}{\min} |e|$. If $m <k$, there exists a $k$-hypergraph $H'= (V',E')$ where $\underset{e \in E}{\min} |e| = m+1$, and having $|E'|\le 2|E|$ and $|V'| \le |V| + 2$ such that Avoider has a winning strategy in Avoider-Enforcer on $H$ if and only if she has one in $H'$.
\end{lemma}

\begin{proof}
Let $H = (V,E)$ be a $k$-hypergraph. Let $m = \underset{e \in E}{\min} |e|$. We define $H' = (V', E')$ as follows. We start from $H' = H$. We add two vertices $\{a_1, a_2\}$ in $V'$. For each edge $e \in E$, we do the following:
\begin{itemize}
    \item If $|e| > m$, we keep $e$ in $E'$.
    \item If $|e| = m$, we consider $e_1 = e \cup \{a_1\}$ and $e_2 = e \cup \{a_2\}$. We then replace $e$ in $E'$ by $e_1$ and $e_2$.
\end{itemize}
We have $|V'|= |V| +2$, $|E'| \le 2|E|$ and $\underset{e \in E'}{\min} |e| = m+1$.

Now, if Avoider (Enforcer resp.) had a winning strategy $\mathcal{S}$ in $E$, we can define a strategy $\mathcal{S}'$ in $E'$ as follows:
\begin{itemize}
    \item If the opponent plays a vertex in $V$, or if it is the first move of the player, play as in $\mathcal{S}$.
    \item If the opponent plays a vertex in $\{a_1,a_2\}$, or if there is no vertex in $V$ available, play an available vertex in $\{a_1,a_2\}$.
\end{itemize}
Following this strategy, Avoider (Enforcer resp.) has played exactly the same vertices in $H'$ as he (she resp.) would have played in $H$ according to $\mathcal{S}$ with the addition of exactly one of $\{a_1,a_2\}$. 

Therefore, if Avoider had a winning strategy in $H$, then for each $e \in E$, there exists one vertex $v \in e$, that Enforcer has played. This vertex is also in $e$ if $e \in E'$ or in $e_1$ and $e_2$ if $e_1$ and $e_2$ were added to $E'$ when we considered $e$. So this strategy ensure Avoider that Enforcer has played one vertex in $e$, and so it is a winning strategy for Avoider.

If Enforcer had a winning strategy in $H$, following this strategy, there exists an edge $e \in E$ in which Avoider has coloured all the vertices. If $|e| \ge m+1$, Avoider has also coloured all the vertices of $e \in E'$, so Enforcer has won. If $|e| = m$, as the strategy $\mathcal{S}'$ forces Avoider to play at least one of $\{a_1, a_2\}$, suppose without loss of generality that she has played $a_1$. Then, she has played $a_1$ and all the vertices of $e$, so she has filled the edge $e_1$. Therefore, this strategy is a winning strategy for Enforcer.

Finally, $H'$ has the same outcome as $H$ and $\underset{e \in E'}{\min} |e| = m+1$.
\end{proof}

\begin{corollary}\label{corollary6-unif}
Avoider-Enforcer is {\sf PSPACE}-complete even on $6$-uniform hypergraphs
\end{corollary}

Several games have been proven to be {\sf PSPACE}-complete in the Maker-Breaker convention, thanks to the proofs of Schaefer or of Rahman and Watson. Due to the similarities between the two convention, some reductions may be adapted to prove that these games are {\sf PSPACE}-complete in the Avoider-Enforcer convention. In particular, we prove in this section that Avoider-Avoider games are {\sf PSPACE}-complete, and we show that the Domination Game and the $H$-Game are {\sf PSPACE}-complete in Avoider-Enforcer convention.

\subsection{Avoider-Avoider games}

We prove here that Avoider-Avoider games are {\sf PSPACE}-complete.

\avoideravoider*

\begin{proof}
    Consider the construction provided in the proof of \cref{corollary6-unif}. Consider $H'$ the hypergraph obtained by adding a vertex $v_0$ in $H$ and adding it in all the hyperedges of $H$. Note that in $H'$, all the hyperedges have size seven. By Lemma~\ref{lemma:included moves}, both player have an optimal strategy in which $v_0$ will be played last, and as the graph has an odd number of vertices, the first player will play it. Therefore, the second player cannot fill a hyperedge and plays as Enforcer would in the Avoider-Enforcer game. By applying the same strategy as in Avoider-Enforcer, If Avoider wins in Avoider-Enforcer, the game ends by a draw, otherwise the second player wins.
\end{proof}

\subsection{ Particular Avoider-Enforcer games}

\subsubsection{Avoider-Enforcer Domination game}

The Maker-Breaker Domination game was introduced by Duchêne {\em et al.} in 2020~\cite{Duchene2020} and follows the study of Domination games on graphs, which were investigated since 2002~\cite{Alon2002, Brevsar2010}. In Maker-Breaker, two players, namely Dominator and Staller alternately claim an unclaimed vertex of the graph, Dominator wins if he manages to take all the vertices from a dominating set. Otherwise, Staller wins. They proved that determining whether Dominator or Staller has a winning strategy is {\sf PSPACE}-complete using a reduction from Maker-Breaker. The Avoider-Enforcer Domination game can be similarly defined, with Staller winning if Dominator claims a dominating set and Dominator winning otherwise. We prove here that determining the winner of the Avoider-Enforcer Domination Game is {\sf PSPACE}-complete. Note that the proof is very similar to the reduction from Maker-Breaker games to Maker-Breaker Domination game.

\begin{problem}{\sc Avoider-Enforcer Domination Game}

Input: A graph $G$

Output: True if and only if Avoider wins the Avoider-Enforcer Domination Game on $G$.
\end{problem}

\begin{restatable}{theorem}{Domination}
Avoider-Enforcer Domination game is {\sf PSPACE}-complete
\end{restatable}

\begin{proof}
The proof of {\sf PSPACE}-completeness in Avoider-Enforcer is roughly the same as in Maker-Breaker.

First, Avoider-Enforcer Domination game is in {\sf PSPACE}, as the number of moves in a game is the number of vertices, and as determining if a set is a dominating set or not, the game is in {\sf PSPACE}.

Let $H = (V_H, E_H)$ be a hypergraph. Without loss of a generality, we can suppose that each vertex is in at least one hyperedge. We construct the following graph $G = (V,E)$ as follows:

\begin{itemize}
    \item For each vertex $u_i$ in $V_H$, we add a vertex $v_i$ in $V$.
    \item For each hyperedge $C$ in $E_H$, we add two vertices $v_C^1$ and $v_C^2$ in $V$.
    \item If a vertex $u_i$ of $V_H$ belongs to a hyperedge $C$ of $E_H$, we add the edges $v_iv_C^1$ and $v_iv_C^2$ to $E$.
\end{itemize}
Note that the graph created here is bipartite

Suppose Avoider (Enforcer resp.) has a winning strategy $\mathcal{S}$ in $H$. We define a strategy $\mathcal{S'}$ for Bob (Alice resp.) in $G$ as follows.

\begin{itemize}
    \item If Avoider (Enforcer resp.) plays the first move in $H$, play first a vertex $v_i$ such that $u_i$ is the first vertex played in $\mathcal{S}$.
    
    Then:
    \item If the opponent plays a vertex $v_i$, play a vertex $v_j$ such that $u_j$ is the answer to the vertex $u_i$ in $\mathcal{S}$.
    \item If a player plays a vertex $v_C^k$ for $k \in \{1,2\}$, plays the vertex $v_C^{k'}$ for $k'\neq k \in \{1,2\}$.
\end{itemize}

Now if Avoider had a winning strategy in $H$, by applying the strategy $\mathcal{S}'$, for any vertex $v_C^i$, Bob has not played all the $v_j$s adjacent to it. Therefore, Alice has played one of them and all the $v_C^i$s are dominated. Moreover, all the vertex $u_j$s are in at least one edge $C$ of $H$. So $v_j$ is dominated by Alice has played one of $(v_C^1, v_C^2)$. So Bob has won.

If Enforcer had a winning strategy in $H$, by applying the strategy $\mathcal{S}'$, Bob knows that there exists a pair of vertices $(v_C^1, v_C^2)$, such that Alice has played all the $v_j$s adjacent to them. As Bob has played exactly one of them, he has not played the second one. Therefore, he does not dominate it and has won. 
\end{proof}

\begin{figure}[ht]
\centering
\scalebox{1}{
\begin{tikzpicture}

\node[noeud] (x1) at (0,0){};
\node[noeud] (x2) at (1.25,0){};
\node[noeud] (x3) at (2.5,0){};
\node[noeud] (x4) at (3.75,0){};

\node[above] at (x1) {$v_1$};
\node[above] at (x2) {$v_2$};
\node[above] at (x3) {$v_3$};
\node[above] at (x4) {$v_4$};

\node[noeud] (c1) at (-0.8,-3){};
\node[noeud] (c1b) at (0.1,-3){};
\node[noeud] (c2) at (1.5,-3){};
\node[noeud] (c2b) at (2.4,-3){};
\node[noeud] (c3) at (3.8,-3){};
\node[noeud] (c3b) at (4.7,-3){};

\node[below=1pt] at (c1){$v_A^1$};
\node[below=1pt] at (c1b) {$v_A^2$};
\node[below=1pt] at (c2) {$v_B^1$};
\node[below=1pt] at (c2b) {$v_B^2$};
\node[below=1pt] at (c3) {$v_C^1$};
\node[below=1pt] at (c3b) {$v_C^2$};

\draw (x1) -- (c1);
\draw (x1) -- (c1b);
\draw (x1) -- (c2);
\draw (x1) -- (c2b);
\draw (x2) -- (c1);
\draw (x2) -- (c1b);
\draw (x2) -- (c3);
\draw (x2) -- (c3b);
\draw (x3) -- (c3);
\draw (x3) -- (c3b);
\draw (x4) -- (c2);
\draw (x4) -- (c2b);
\draw (x4) -- (c3);
\draw (x4) -- (c3b);

\node[noeud] (u1) at (-5.5,-2.25){};
\node[noeud] (u2) at (-5.5,-0.75){};
\node[noeud] (u3) at (-4,-0.75){};
\node[noeud] (u4) at (-4,-2.25){};
\node (u1') at (-5.3,-2.25){};
\node (u2') at (-5.3,-0.75){};
\node (u3') at (-3.8,-0.75){};
\node (u4') at (-3.8,-2.25){};

\node[right] at (u1){$u_1$};
\node[right] at (u2){$u_2$};
\node[right] at (u3){$u_3$};
\node[right] at (u4){$u_4$};

  \draw \convexpath{u2',u3',u4'}{0.5cm};
  \draw \convexpath{u1',u2'}{0.45cm};
  \draw \convexpath{u1',u4'}{0.4cm};

\node[left] at (-5.7, -1.4){$A$};
\node[right] at (-3.4, -1){$B$};
\node at (-4.5, -2.9){$C$};

\end{tikzpicture}
}
\caption{Reduction from {\sc Avoider-Enforcer}.}
\label{fig:Domination}
\end{figure}
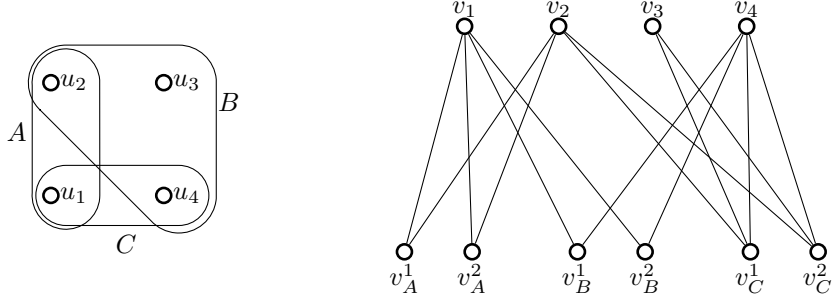

In Figure~\ref{fig:Domination} we provide an example of reduction. Note that by connecting all the vertices $v_i$, they will form a clique, and this gives the proof Avoider-Enforcer Domination game on split graphs. So the Avoider-Enforcer game is {\sf PSPACE}-complete on bipartite and split graphs, as the Maker-Breaker Domination game for roughly the same construction.

\subsubsection{Vertex H-game}

The vertex $H$-Game has been introduced by Kronenberg, Mond and Naor in \cite{Kronenberg2019} on random graphs. It is presented in several conventions, but we will focus here on the Avoider-Enforcer one. The game is played as follows:

Let $H$ be a graph. Avoider and Enforcer play on the vertex set of another graph $G$. Alternately, Avoider and Enforcer claim an unclaimed vertex of $G$. Avoider wins if the set of vertices she has claimed do not contain $H$ as a subgraph (not necessarily induced). Otherwise, Enforcer wins. 

We prove here that determining the winner of the Avoider-Enforcer $H$-Game is a {\sf PSPACE}-complete problem for several graphs $H$ .

\begin{problem}{\sc Avoider-Enforcer $H$-Game}

\noindent Input: A graph $G$

\noindent Output: True if and only if Avoider wins the Avoider-Enforcer $H$-Game played on $G$.
\end{problem}

We prove that the vertex $H$-game is {\sf PSPACE}-complete  for several graphs $H$. To define these graphs, we first need to define some graphs and operations:

\begin{itemize}
    \item We will denote by $I_k$ the graph being an independent set of size $k$, i.e. containing $k$ vertices and no edge.
    \item If $G$ and $H$ are two graphs, we denote by $G\bowtie H$ their join, i.e. if $G = (V_G, E_G)$ and $H = (V_H, E_H)$, we have $G \bowtie H = (V,E)$ with $V = V_G \cup V_H$ and $E = E_G \cup E_H \cup \{(v_G,v_H)| v_G\in V_G, v_H \in V_H\}$.
    \item If $G$ and $H$ are two graphs, we denote by $G\boxtimes H$ their strong product, i.e. if $G = (V_G, E_G)$ and $H = (V_H, E_H)$, we have $G \boxtimes H = (V,E)$ with $V = \{ x_{u,v}| u \in V_G, v \in V_H \}$ and $E = \{ (x_{u_1,v_1}, x_{u_2,v_2}) | \left (u_1 = u_2 \mbox{ or } (u_1,u_2) \in E_G \right )$ and $\left ( v_1=v_2 \mbox{ or } (v_1,v_2) \in E_H \right ) \}$.
\end{itemize}

Remark that for any graph $G$, $G \boxtimes P_2$ (where $P_2$ design the path of length $2$) is obtained by taking two copies of $G$ and connecting each vertex to its copy and its copy's neighbours.

We prove here that determining the winner of the Avoider-Enforcer $H$-Game is a {\sf PSPACE}-complete problem for several graphs $H$ .

\begin{theorem}
Let $H_0$ be a graph containing at least one edge or at least $6$ vertices, and let $k\ge 6$. Consider $H = I_k \bowtie H_0 $. The {\sc Avoider-Enforcer $H$-game} problem is {\sf PSPACE}-complete.  
\end{theorem}

Note that complete bipartite graphs $K_{n,m}$, with $n,m \ge 6$, are of this type. Indeed, $K_{n,m} = I_k \bowtie H_0$ for $H_0 = I_m$ and $k= n$.

\begin{proof}
First, the $H$-game is in {\sf PSPACE}. Indeed, as it is a positional game, if $G = (V,E)$ is a graph, the game ends after at most $|V|$ moves. After that, determining whether a graph $H$ is a subgraph of a graph $G$ can be done in polynomial space.

We do our reduction from {\sc Avoider-Enforcer} on $6$-uniform hypergraphs which is {\sf PSPACE}-complete by \cref{corollary6-unif}. Let $H_0$ be a graph containing at least one edge or at least $6$ vertices, and let $k \ge 6$. Let $H = I_k \bowtie H_0$. To avoid confusion while describing the strategies in the two games, we will call the players of the $H$-Game Alice and Bob, with Alice avoiding creating a subgraph $H$ and Bob forcing her to create one.

Let $H' = (V', E')$ be a $6$-uniform hypergraph. Let $H_0'$ be the strong product $H_0 \boxtimes P_2$.
We build $G = (V,E)$ an instance of $H$-Game as follows: 

\begin{itemize}
    \item {\bf Step 1}: For any vertex $v_i' \in V'$, we add a vertex $v_i \in V$. 
    \item {\bf Step 2}: For any edge $C \in E'$, we add $2(k-6)$ vertices $v^C_1, \dots, v^C_{2(k-6)}$ (note that if $k=6$ these vertices do not exist).
    \item {\bf Step 3}: For any edge $C \in E'$, we add a copy $H^C_0$ of the graph $H_0'$ in $G$, and we connect any vertex of $H^C_0$ to all the vertices $v_i$ such that $v_i' \in C$ and to all the vertices $v^C_j$ for $1 \le j \le 2k-6$.
\end{itemize}

\begin{claim}
If Avoider has a winning strategy in $H'$, Alice has a winning strategy in $G$.
\end{claim}

\begin{claimproof}
Suppose Avoider has a winning strategy $\mathcal{S}'$ in $H'$. We define Alice's strategy $\mathcal{S}$ in $G$ as follows:

\begin{itemize}
    \item She starts by playing the vertex $v_i$, corresponding to the vertex $v'_i$ that Avoider would have played in $H'$ according to $\mathcal{S}'$.
    \item If Bob plays a vertex $v_i$, she answers with the vertex $v_{j}$ corresponding to the vertex $v_j'$ that Avoider would have played by $\mathcal{S}'$ in $H'$ if Enforcer has played $v_i'$.
    \item In $H^C_0$, as it is a strong product $H_0 \boxtimes P_2$, Alice considers the pairing between any vertex and its copy in the strong product. If Bob plays one of them, she plays the second one.
    \item For any edge $C$ in $H'$, Alice considers the set of vertices $v^C_j$ for $1 \le j \le 2k-6$. As there is an even number of them, if Bob plays one of them, she can also play one of them.
\end{itemize}

At the end of the game, by the matching strategy, for each clause $C$ in $H'$, Alice will have played  exactly a copy of $H_0$ in each $H^C_0$, exactly $k-6$ vertices among the $v^C_j$ with $1 \le j \le 2(k-6)$, and the vertices $v_i$ corresponding to the $v_i'$ that Avoider would have played according to $\mathcal{S}'$ in $H'$.

Now, consider any copy $H_1$ of $H$ in $G$. Suppose that Alice has played all the vertices of $H_1$.

Suppose $H_0$ has at least one edge. We first prove that $H_1$ cannot contain two vertices $v_C$ and $v_{C'}$ for $C \neq C'$ that have been created by the Steps~$2$ or Step~$3$ of our construction. Suppose it does. Consider a decomposition of $H_1 = I_k \bowtie H^1_0$. By construction, the $v_C$s and the $v_{C'}$s are not adjacent. Therefore, as these two components are fully connected one to the other, they must either be both in $I_k$ or both in $H^1_0$.

\noindent{\bf Case 1:} $H_0$ has at least one edge.

$H_1$ also has one edge $e = (u_1,u_2)$. As $I_k$ is a stable set, we have $u_1,u_2 \in H^1_0$. Now, as $u_1$ can only be adjacent to vertices $v_i$s created by the vertices of $H'$(Step~1), $v_C$ and $v_{C}'$ cannot be both adjacent to $u_1$. Therefore, as we supposed $C \neq C'$, and as only vertices created during Step~1 can be adjacent to both $v_C$ and $v_{C'}$, any vertex in $I_k$ must be a vertex $v_i$. Which is not possible otherwise, Alice would have play $k \ge 6$ vertices $v_i$ adjacent to a same $v_C$, which means that, according to $\mathcal{S}'$ she would have played $k \ge 6$ vertices in the same hyperedge in $H'$, which contradicts the fact that $\mathcal{S}'$ was a winning strategy for Avoider in $H'$.

Now, as all the vertices of $H_1$ are either $v_i$s or were created by considering the same hyperedge $C$, and $|H_1| = |H_0| + k $, by the pairing, we know that exactly $6$ of them are $v_i$s. As there are no edges between the $v_i$s, they must all be on the same side of the join, and therefore, they are all connected to a same vertex. By construction, this is only the case if these six vertices were in a same hyperedge of $H'$ which contradicts that $\mathcal{S}'$ was a winning strategy for Avoider in $H'$.

\noindent{\bf Case 2:} $H_0$ has no edges and has $k'\ge 6$ vertices.

This means that $H$ can be written $I_k \bowtie I_{k'}$ for $k,k' \ge 6$ (note that $H$ is a complete bipartite graph). Once again, consider two vertices $v_C$ and $v_{C'}$ for $C \neq C'$ in $H_1$. As they cannot be adjacent, they must be both in $I_k$ or both in $I_{k'}$. Thus, $v_C$ and $v_{C'}$ have $\min(k,k') \ge 6$ common neighbours. This implies that, if $C \neq C'$, at least one of their common neighbour is not a vertex $v_i$ created during Step~1, otherwise Alice would have played six vertices in the same hyperedge of $H'$. This is not possible by construction. So once again, $H_1$ cannot contain $v_C$ and $v_{C'}$ created from different hyperedges from $H'$. Now, if Alice has played all the vertices of $H_1$, by construction as $|H_1| = k+k'$, and as her pairing strategy ensure her to play $k-6$ vertices created during step~2 and $k'$ created during step~3, necessarily, she has played six vertices $v_i$ creating during step~1. As there are no edges between these six vertices, they must all be in the same independent set $I_k$ or $I_{k'}$. Thus, they have a common neighbour. This common neighbour must then be a vertex $v_C$ creating during step~3 as only them are connected to the $v_i$s. Finally, these six vertices corresponds to six vertices $v_i'$s that are in the same hyperedge $C$ of $H'$. Once again, this contradicts the fact that $\mathcal{S}'$ was a winning strategy for Avoider in $H'$.

\end{claimproof}

\begin{claim}
If Enforcer has a winning strategy $\mathcal{S}'$ in $H'$, Bob has a winning strategy in $G$.
\end{claim}

\begin{claimproof}

Let $\mathcal{S}'$ be a winning strategy for Enforcer in $H'$. We consider a strategy $\mathcal{S}$ for Bob in $G$ as follows:

\begin{itemize}
    \item If Alice plays a vertex $v_i$, he answers with the vertex $v_j$ that corresponds to the vertex $v_j'$ that Enforcer would have played in response to $v_i'$ in $\mathcal{S}'$.
    \item In $H^C_0$, as it is a strong product $H_0 \boxtimes P_2$, Bob considers the pairing between any vertex and its copy in the strong product. If Alice plays one of them, he plays the second one.
    \item For any edge $C \in H'$, Bob pairs the vertices $v^C_j$ for $1 \le j \le 2k-6$. If Alice plays one of them, he plays one of them too.
    \item If at a certain moment of the game, it is Bob's turn and the remaining vertices are all in some $H^C_0$ or vertices $v^C_j$s, he applies the pairing strategy, so that Alice plays once in any pair of the matching by Lemma~\ref{pairing strategy}.
\end{itemize}

Consider the graph at the end of the game. As $\mathcal{S}$ was a winning strategy for Enforcer in $H'$, there exists a hyperedge $C \in H'$ in which Avoider has played the six vertices. Up to a renaming of the vertices, denote by $v_1', \dots, v_6'$ be these six vertices. Alice has then played $v_1, \dots, v_6$ in $G$. According to the pairing strategy, Bob knows that Alice will play exactly $k-6$ vertices from the $v^C_j$, denote them $v_7, \dots, v_k$, and exactly one copy $H_1$ of $H_0$ from the vertices of $H^C_0$. Now, by construction, the vertices $v_1, \dots v_k$ are a stable set and all the edges exist between any $v_i$ ($1 \le i \le k$) and any vertex $v$ of $H_1$. Thus, the subgraph formed by these vertices, which were all played by Alice, is isomorphic to $I_k \bowtie H_0 = H$. Thus, Bob has won.
\end{claimproof}

Finally, the $H$-Game played on $G$ is won by Alice if and only if the Avoider-Enforcer game played on $H'$ is won by Avoider, and determining the winner of the $H$-Game is {\sf PSPACE}-complete. 
\end{proof}

Theses two games are two examples of games in which the proof of {\sf PSPACE}-completeness in Avoider-Enforcer convention is similar to the one in Maker-Breaker convention, in several other games the Maker-Breaker complexity proof can be adapted to an Avoider-Enforcer one.

\section*{Acknowledgments}
We want to thank Eric Duchêne and Aline Parreau for their help in the redaction of this paper

\bibliography{main}

\end{document}